\documentclass[a4paper,11pt]{article}
\usepackage{amsmath,amsfonts,amsthm,dsfont,amssymb}
\usepackage[blocks]{authblk}
\usepackage{cite}
\usepackage{graphicx}
\usepackage{amssymb,array,arydshln}
\usepackage{amsfonts}
\usepackage{amsbsy}
\usepackage{cite}
\usepackage{amssymb}
\usepackage{amsmath}
\usepackage{graphics}
\newenvironment{keywords}{\noindent\textbf{Keywords:}}{}
\newenvironment{classification}{\noindent\textbf{AMS subject classifications.}}{}
\date{}
\newcommand{\email}[1]{\texttt{\small #1}}
\newtheorem{theorem}{Theorem}[section]
\newtheorem{remark}[theorem]{Remark}
\newtheorem{example}[theorem]{Example}
\newtheorem{lemma}[theorem]{Lemma}
\newtheorem{corollary}[theorem]{Corollary}
\newtheorem{definition}[theorem]{Definition}
\newtheorem{proposition}[theorem]{Proposition}

\newtheorem{conjecture}{Conjecture}
\newtheorem{question}{Question}
\newcommand{\bea}{\begin{eqnarray}}
\newcommand{\eea}{\end{eqnarray}}
\newcommand{\beq}{\begin{eqnarray*}}
\newcommand{\eeq}{\end{eqnarray*}}
\def \bd{\begin{definition}}
\def \ed{\end{definition}}
\def \bqu{\begin{question}}
\def \equ{\end{question}}
\def \bcc{\begin{conjecture}}
\def \ecc{\end{conjecture}}
\def \bt{\begin{theorem}}
\def \et{\end{theorem}}
\def \bl{\begin{lemma}}
\def \el{\end{lemma}}
\def \bc{\begin{corollary}}
\def \ec{\end{corollary}}
\def \be{\begin{equation}}
\def \ee{\end{equation}}
\def \ben{\begin{enumerate}}
\def \een{\end{enumerate}}
\def \ba{\begin{array}}
\def \ea{\end{array}}
\def \bp{\begin{proposition}}
\def \ep{\end{proposition}}
\def \bx{\begin{example}}
\def \ex{\end{example}}
\def \br{\begin{remark}}
\def \er{\end{remark}}
\def \bdsc{\begin{description}}
\def \edsc{\end{description}}

\def\1{1\!\!1}
\def\0{0\!\!0}

\begin{document}
	\title{A note on the distance spectra of co-centralizer graphs
	}

\author[1]{Jharna Kalita}
\author[2]{Somnath Paul\footnote{Corresponding  Author.}}
\affil[1]{Department of Applied Sciences\\ Tezpur University\\ Napaam-784028, Assam, India. \email{app21104@tezu.ac.in}}
\affil[2]{Department of Applied Sciences\\ Tezpur University\\ Napaam-784028, Assam, India. \email{som@tezu.ernet.in}}

\pagestyle{myheadings} \markboth{J. Kalita \& S. Paul}{A note on the distance spectra of co-centralizer graphs}
   \maketitle
\begin{abstract}
Let $G$ be a finite non abelian group. The centralizer graph of $G$  is a simple undirected graph $\Gamma_{cent}(G)$, whose vertex set consists of proper centralizers of $G$ and two vertices are adjacent if and only if their cardinalities are identical {\rm\cite{omer}}. We call the complement of the centralizer graph as the co-centralizer graph. In this paper, we investigate the distance, distance (signless) Laplacian spectra of co-centralizer graphs of some classes of finite non-abelian groups, and obtain some conditions on a group so that the co-centralizer graph is distance, distance (signless) Laplacian integral.
\end{abstract}


\section{Introduction}

Let $G$ be a finite non-ableian group. The centralizer graph of $G$  introduced in \cite{omer}, is a simple undirected graph $\Gamma_{cent}(G)$, whose vertex set consists of proper centralizers of $G$ and two vertices are adjacent if and only if their cardinalities are identical. The authors in \cite{omer} have discussed the structure of the centralizer graph for Dihedral group, quaternion group and dicyclic group.

A different definition of centralizer graph is given in \cite{gazvin}, where two vertices $x$ and $y$ of $G\setminus Z(G)$ are adjacent if their centralizers are equal. They have discussed some graph theoretical properties there. We define the complement of the centralizer graph as \textit{the co-centralizer graph} and denote it  by $\overline{\Gamma_{cent}(G)}$ .

For a simple graph $H,$ the \textit{distance} between two of its vertices $u$ and $v$ is defined as the length of a shortest path between $u$ and $v$ in it, and is denoted by $d_{uv}.$ The \textit{distance matrix} of $H$ is denoted by $D(H),$ with the $(u,v)$-entry equal to $d_{uv}.$

The \textit{transmission} of a vertex $v$ is defined to be the sum of the distances from $v$ to all other vertices in the graph. Let $Tr(H)$ denote the diagonal matrix with $i$-th diagonal entry being the  transmission of the $i$-th vertex in $H.$ In \cite{AH2}, M. Aouchiche and P. Hansen have defined the Laplacian and the signless Laplacian for the distance matrix of a connected graph, which are analogous to the Laplacian matrix and signless Laplacian matrix for adjacency matrix. The matrix $ D^L(H)=Tr(H)-D(H) $ is called the \textit{distance Laplacian matrix} of $H$, while the matrix  $D^{Q}(H)=Tr(H)+D(H)$ is called the \textit{distance signless Laplacian matrix of $H$}. For a beautiful survey of distance matrix, the reader can see \cite{AH1}.

If $M$ is a symmetric matrix, then the characteristic polynomial of $M$ has only real zeroes. We will represent this family of eigenvalues (known as the \textit{spectrum}) as $$\sigma_M=\left(
                  \begin{array}{cccc}
                    \mu_1 & \mu_2 &\cdots&\mu_p\\
                    m_1 & m_2 &\cdots & m_p
                  \end{array}
                \right),
$$ where $\mu_1,\mu_2,\ldots,\mu_p$ are the distinct eigenvalues of $M$ and $m_1, m_2,\ldots,m_p$ are the corresponding multiplicities. Since each of $D(H),~D^L(H)$ and $D^Q(H)$ is symmetric, we will refer the corresponding spectrum as the \textit{distance, distance Laplacian} and \textit{distance signless Laplacian spectrum}, respectively.

A graph is called distance (respectively distance (signless) Laplacian) integral if the distance (respectively distance (signless) Laplacian) spectrum consists entirely of integers.  In this article, we consider some finite non-abelian groups, namely the generalized quaternion group, the dihedral group, the quasidihedral group, the metacyclic group, and the projective special linear group,  and investigate the distance, distance (signless) Laplacian spectra of co-centralizer graphs of them. Moreover, we obtain some conditions so that their co-centralizer graph is distance, distance (signless) Laplacian integral.

\section{Preliminaries}

Consider $Q_{4n} =<x,y : x^{2n} = 1, x^n = y^2, yx=x^{-1}y>,$ be the generalized quaternion group of order $4n,$  where $n\ge2$ , and $Z(Q_{4n})$= $\{1,x^n\}$ be its centre. Then $Q_{4n}$ can be written as $A\cup B,$ where $A=\{1,x,x^2,\ldots,x^{2n-1}\}$ and $B=\{y,xy,x^2y,\ldots,x^{2n-1}y\},$ where each element of $B$ is of order 4. We note that for any $z\in Z(Q_{4n}),$ and $1\le i\le 2n-1,$
\begin{eqnarray*}
  C_{Q_{4n}}(x)
  &=&C_{Q_{4n}}(x^iz)\\ &=& Z(Q_{4n})\cup x Z(Q_{4n})\cup x^2Z(Q_{4n})\cup\ldots\cup x^{n-1}Z(Q_{4n})\\
    &=& \{1,x^n\} \cup x \{1,x^n\}\cup x^2\{1,x^n\}\cup\ldots\cup x^{n-1}\{1,x^n\}\\
    &=& \{1,x^n\} \cup \{x,x^{n+1}\}\cup \{x^2,x^{n+2}\}\cup\ldots\cup \{x^{n-1},x^{2n-1}\} \\
    &=& \{1,x,x^2,\ldots,x^{2n-1}\}.
\end{eqnarray*}
Moreover, for $1\le j\le n,$
\begin{eqnarray*}
  C_{Q_{4n}}(yx^j)=C_{Q_{4n}}(yx^jz) &=& Z(Q_{4n})\cup yx^j Z(Q_{4n}) \\
    &=& \{1,x^n\} \cup yx^j \{1,x^n\}\\
    &=& \{1,x^n\} \cup \{yx^j,yx^{n+j}\}.
\end{eqnarray*}

Therefore, $Q_{4n}$ has $n+1$ distinct centralizers. One of them has cardinality $2n$ and others have cardinality $4$. Thus, from the definition of centralizer graph, it follows that $\Gamma_{cent}(Q_{4n})$ is a graph with $n+1$ vertices where one component is $K_n$ and another is an isolated vertex. Therefore the co-centralizer graph of $Q_{4n}$ is $K_{1,n}$.

The following result gives the distance characteristics polynomial for $K_{n_1,n_2,\ldots,n_k}$ and will be useful to derive some of our main results.
\bl{\rm \cite{disc}}\label{lem1} The distance characteristics polynomial of the complete multipartite graph $K_{n_1,n_2,\ldots,n_k}$ is
\begin{equation}\label{e1}
 P_D(\lambda)=(\lambda+2)^{n-k}\left[\prod_{i=1}^k(\lambda-n_i+2)-\sum_{i=1}^k n_i\prod_{j=1,j\neq i}^k(\lambda-n_j+2)\right].
\end{equation} \el
\section{Spectra of $\overline{\Gamma_{cent}(Q_{4n})}$}\label{sec1}

In this section, we consider the co-centralizer graph of $Q_{4n}$ and obtain the distance, distance Laplacian and distance signless Laplacian spectra of it. As it is already observed in the previous section,  $\overline{\Gamma_{cent}(Q_{4n})}=K_{1,n}$.
\subsection{Distance spectrum of $\overline{\Gamma_{cent}(Q_{4n})}$}
By equation (\ref{e1}) we get,
\begin{eqnarray*}
  &&P_{D(\overline{\Gamma_{cent}(Q_{4n})})}(\lambda) \\
  &=&(\lambda+2)^{(n+1)-2}[(\lambda-1+2)(\lambda-n+2)-n(\lambda-1+2)-(\lambda-n+2)]\\
  &=&(\lambda+2)^{n-1}[(\lambda+1)(\lambda-n+2)-n(\lambda+1)-(\lambda-n+2)]\\
  &=&(\lambda+2)^{n-1}[\lambda^2+(2-2n)\lambda-n]
\end{eqnarray*}

Therefore, eigenvalues of $D(\overline{\Gamma_{cent}(Q_{4n})})$ are $-2$ with multiplicity $(n-1),$ and two roots of the equation $\lambda^2+(2-2n)\lambda-n=0$, ie $(n-1)\pm\sqrt{n^2-n+1}.$ Hence, we have the following theorem.

\bt\label{dt1} Let $\overline{\Gamma_{cent}(Q_{4n})}$ be the co-centralizer graph of the generalized quaternion group of order $4n.$ Then
\begin{enumerate}
\item[{\rm(a)}] $-2\in\sigma(D(\overline{\Gamma_{cent}(Q_{4n})}))$ with multiplicity $n-1;$
\item[{\rm(b)}] $(n-1)\pm\sqrt{n^2-n+1}\in\sigma(D(\overline{\Gamma_{cent}(Q_{4n})})),$ each with multiplicity 1.
\end{enumerate}\et
Therefore, by Theorem~\ref{dt1}, it follows that $\overline{\Gamma_{cent}(Q_{4n})}$ is distance integral if  $n^2-n+1$ is a perfect square.

\subsection{Distance Laplacian and distance signless Laplacian spectrum of $\overline{\Gamma_{cent}(Q_{4n})}$}\label{subi}
It is already observed in {\rm\cite{AH2}} that, the distance Laplacian spectum and the distance signless Laplacian spectrum of a star graph $K_{1,n-1}$ is: $\sigma_{D^L(K_{1,n-1})}=\left(
                                                                                                                   \begin{array}{ccc}
                                                                                                                     2n-1 & n & 0 \\
                                                                                                                     n-2 & 1 & 1 \\
                                                                                                                   \end{array}
                                                                                                                 \right)
,$ and $\sigma_{D^Q(K_{1,n-1})}=\left(
                                                                                                                   \begin{array}{cc}
                                                                                                                     2n-5 & \frac{\left[(5n-8)\pm\sqrt{9n^2-32n+32}\right]}{2} \\
                                                                                                                     n-2 & 1 \\
                                                                                                                   \end{array}
                                                                                                                 \right).$
Therefore the distance Laplacian spectrum of $\overline{\Gamma_{cent}(Q_{4n})}$ is $\left(
                                                                                                                   \begin{array}{ccc}
                                                                                                                     2n+1 & n+1 & 0 \\
                                                                                                                     n-1 & 1 & 1 \\
                                                                                                                   \end{array}
                                                                                                                 \right)
.$
And the distance signless Laplacian spectrum of $\overline{\Gamma_{cent}(Q_{4n})}$ is: $\left(
                                                                                                                   \begin{array}{cc}
                                                                                                                     2n-3 & \frac{\left[(5n-3)\pm\sqrt{9n^2-14n+9}\right]}{2} \\
                                                                                                                     n-1 & 1 \\
                                                                                                                   \end{array}
                                                                                                                 \right).$

Therefore, $\overline{\Gamma_{cent}(Q_{4n})}$ is always distance Laplacian integral but it is distance signless Laplacian integral only if $(5n-3)\pm\sqrt{9n^2-14n+9}$ is even.

\section{Spectra of $\overline{\Gamma_{cent}(D_{2m})}$}
In this section, we consider the co-centralizer graph of the dihedral group $D_{2m}=<a,b:a^m=b^2=1,bab^{-1}=a^{-1}>,$  and obtain the distance, distance Laplacian and distance signless Laplacian spectra of it. It follows from {\rm\cite{rkn}} that if $m$ is even, then $D_{2m}$ has $\frac{m}{2}+1$ distinct centralizers with one of cardinality $m$ and the others are of cardinality $4$. Therefore the co-centralizer graph of $D_{2m}$ is $\overline{\Gamma_{cent}(D_{2m})}= K_{1,\frac{m}{2}}.$ Similarly, if $m$ is odd, then $D_{2m}$ has $m+1$ distinct centralizers with one of cardinality $m$ and the rests are of cardinality $2$. Therefore, $\overline{\Gamma_{cent}(D_{2m})}=K_{1,m}$. Hence, similar to Theorem~\ref{dt1}, we have the following.
\bt\label{ddt1a} Let $\overline{\Gamma_{cent}(D_{2m})}$ be the non centralizer graph of the dihedral group of order $2m$. If $m$ is even, then  \begin{enumerate}
\item[{\rm(a)}] $-2\in\sigma(D(\overline{\Gamma_{cent}(D_{2m})}))$ with multiplicity $\frac{m}{2}-1;$
\item[{\rm(b)}] $(\frac{m}{2}-1)\pm \frac{1}{2}\sqrt{m^2-2m+4} \in\sigma(D(\overline{\Gamma_{cent}(D_{2m})}),$ each with multiplicity 1.
\end{enumerate}
And if $m$ is odd,then
\begin{enumerate}
\item[{\rm(a)}] $-2\in\sigma(D(\overline{\Gamma_{cent}(D_{2m})}))$ with multiplicity $m-1;$
\item[{\rm(b)}] $(m-1)\pm\sqrt{m^2-m+1} \in\sigma(D(\overline{\Gamma_{cent}(D_{2m})})),$ each with multiplicity 1.
\end{enumerate}
 \et
Therefore, if $m$ is even, then by Theorem~\ref{ddt1a}, $\overline{\Gamma_{cent}(D_{2m})}$ is distance integral if  $\sqrt{m^2-2m+4} $ is even. Also if $m$ is odd, then $\overline{\Gamma_{cent}(D_{2m})}$ is distance integral if $m^2-m+1$ is perfect square.

Moreover, as discussed in Subsection~\ref{subi}, the following two results holds.
\bt\label{dddg2t1} Let $\overline{\Gamma_{cent}(D_{2m})}$ be the co-centralizer graph of the dihedral group of order $2m$. If $m$ is even, then
\begin{enumerate}
\item[{\rm(a)}] $0\in\sigma(D^L(\overline{\Gamma_{cent}(D_{2m}))})$ with multiplicity $1;$
\item[{\rm(b)}] $(\frac{m}{2}+1)\in\sigma(D^L(\overline{\Gamma_{cent}(D_{2m})})),$ with multiplicity $1;$
\item[{\rm(c)}] $(m+1)\in\sigma(D^L(\overline{\Gamma_{cent}(D_{2m})})),$ with multiplicity $\frac{m}{2}-1;$

\end{enumerate}
And if $m$ is odd, then
\begin{enumerate}
\item[{\rm(a)}] $0\in\sigma(D^L(\overline{\Gamma_{cent}(D_{2m}))})$ with multiplicity $1;$
\item[{\rm(b)}] $(m+1)\in\sigma(D^L(\overline{\Gamma_{cent}(D_{2m})})),$ with multiplicity $1.$
\item[{\rm(c)}] $(2m+1)\in\sigma(D^L(\overline{\Gamma_{cent}(D_{2m})})),$ with multiplicity $m-1.$\end{enumerate}\et

\bt\label{tdg2t1} Let $\overline{\Gamma_{cent}(D_{2m})}$ be the co-centralizer graph of the dihedral group of order $2m$. If $m$ is even, then
\begin{enumerate}
\item[{\rm(a)}] $m-3\in\sigma(D^Q(\overline{\Gamma_{cent}(D_{2m}))})$ with multiplicity $\frac{m}{2}-1;$
\item[{\rm(b)}] $\frac{1}{2}\left[(\frac{5m}{2}-3)\pm\sqrt{\frac{9m^2}{4}-7m+9}\right]\in\sigma(D^Q(\overline{\Gamma_{cent}(D_{2m})})),$ each with multiplicity $1;$

\end{enumerate}
And, if $m$ is odd,
\begin{enumerate}
\item[{\rm(a)}] $2m-3\in\sigma(D^Q(\overline{\Gamma_{cent}(D_{2m}))})$ with multiplicity $m-1;$
\item[{\rm(b)}] $\frac{1}{2}\left[(5m-3)\pm\sqrt{9m^2-14m+9}\right]\in\sigma(D^Q(\overline{\Gamma_{cent}(D_{2m})})),$ each with multiplicity $1.$
\end{enumerate}\et
Thus, $\overline{\Gamma_{cent}(D_{2m}))}$ is always distance Laplacian integral. Also  if $m$ is even, then by Theorem~\ref{tdg2t1}, $\overline{\Gamma_{cent}(D_{2m})}$ is distance signless Laplacian integral if  $(\frac{5m}{2}-3)\pm\sqrt{\frac{9m^2}{4}-7m+9}$ is even. And if $m$ is odd, then $\overline{\Gamma_{cent}(D_{2m})}$ is distance signless Laplacian integral if $(5m-3)\pm\sqrt{9m^2-14m+9}$ is even.

\section{Spectra of $\overline{\Gamma_{cent}(QD_{2^n})}$}

In this section, we consider the co-centralizer graph of the Quasidihedral group $QD_{2^n}=<a,b:a^{2^{n-1}}=b^2=1,bab^{-1}=a^{2^{n-2}-1}>,$ where $n\ge4,$ and obtain the distance, distance Laplacian and distance signless Laplacian spectra of it. It has $2^{n-2}+1$ distinct centralizers, of which one has cardinality $2^{n-1}$ and the rests are of $4$. Therefore the co-centralizer graph of $QD_{2^{n}}$ is: $\overline{\Gamma_{cent}(QD_{2^n})}= K_{1,2^{n-2}}$. Thus, similar to Theorem~\ref{dt1} we have the following.

\bt\label{dg2t1} Let $\overline{\Gamma_{cent}(QD_{2^n})}$ be the co-centralizer graph of the Quasidihedral group of order $2^n.$ Then
\begin{enumerate}
\item[{\rm(a)}] $-2\in\sigma(D(\overline{\Gamma_{cent}(QD_{2^n}))})$ with multiplicity $2^{n-2}-1;$
\item[{\rm(b)}] $(2^{n-2}-1)\pm\sqrt{2^{2n-4}-2^{n-2}+1}\in\sigma(D(\overline{\Gamma_{cent}(QD_{2^n})})),$ each with multiplicity 1.
\end{enumerate}\et

Therefore, by Theorem~\ref{dg2t1}, it follows that $\overline{\Gamma_{cent}(QD_{2^n})}$ is distance integral if  $2^{2n-4}-2^{n-2}+1$ is a perfect square. Similar to Theorem~\ref{dddg2t1} and Theorem~\ref{tdg2t1}, we have the following two results.

\bt\label{dg2t2} Let $\overline{\Gamma_{cent}(QD_{2^n})}$ be the co-centralizer graph of the Quasidihedral group of order $2^n.$ Then
\begin{enumerate}
\item[{\rm(a)}] $0\in\sigma(D^L(\overline{\Gamma_{cent}(QD_{2^n}))})$ with multiplicity $1;$
\item[{\rm(b)}] $(2^{n-2}+1)\in\sigma(D^L(\overline{\Gamma_{cent}(QD_{2^n})})),$ with multiplicity $1;$
\item[{\rm(c)}] $(2^{n-1}+1)\in\sigma(D^L(\overline{\Gamma_{cent}(QD_{2^n})})),$ with multiplicity $2^{n-2}-1.$
\end{enumerate}\et

\bt\label{dg2t3} Let $\overline{\Gamma_{cent}(QD_{2^n})}$ be the co-centralizer graph of the Quasidihedral group of order $2^n.$ Then
\begin{enumerate}
\item[{\rm(a)}] $(2^{n-1}-3)\in\sigma(D^Q(\overline{\Gamma_{cent}(QD_{2^n}))})$ with multiplicity $2^{n-2}-1;$
\item[{\rm(b)}] $\frac{1}{2}\left[(5(2^{n-2})-3)\pm\sqrt{9(2^{2n-4})-14(2^{n-2})+9} \right]\in\sigma(D^Q(\overline{\Gamma_{cent}(QD_{2^n})})),$ each with multiplicity $1;$
\end{enumerate}\et

This shows that $\overline{\Gamma_{cent}(QD_{2^n})}$ is always distance Laplacian integral but it is distance signless Laplacian integral only if $\left(5(2^{n-2})-3\right)\pm\sqrt{9(2^{2n-4})-14(2^{n-2})+9}$ is even.

\section{Spectra of $\overline{\Gamma_{cent}(M_{2mn})}$}

In this section, we consider the co-centralizer graph of the  Metacyclic group $M_{2mn}=<a,b: a^m=b^{2n}=1,bab^{-1}=a^{-1}>,$ and obtain the distance, distance Laplacian and distance signless Laplacian spectra of it, where $m>2.$ It can be easily observed that  $\overline{\Gamma_{cent}(M_{2mn})}\cong\left\{
                                         \begin{array}{ll}
                                           K_{1,m}, & \hbox{if $m$ is odd;} \\
                                           K_{1,\frac{m}{2}}, & \hbox{if $m$ is even.}
                                         \end{array}
                                       \right.$
Thus, similar to Theorem~\ref{dt1}, we have the following.
\bt\label{dm2t1a} Let $\overline{\Gamma_{cent}(M_{2mn})}$ be the co-centralizer graph of the Metacyclic group of order $2mn$. If $m$ is even, then
\begin{enumerate}
\item[{\rm(a)}] $-2\in\sigma(D(\overline{\Gamma_{cent}(M_{2mn}))})$ with multiplicity $\frac{m}{2}-1;$
\item[{\rm(b)}] $(\frac{m}{2}-1)\pm\frac{1}{2}\sqrt{m^2-2m+4}\in\sigma(D(\overline{\Gamma_{cent}(M_{2mn})})),$ each with multiplicity 1.
\end{enumerate}
And if $m$ is odd, then
\begin{enumerate}
\item[{\rm(a)}] $-2\in\sigma(D(\overline{\Gamma_{cent}(M_{2mn}))})$ with multiplicity $m-1;$
\item[{\rm(b)}] $(m-1)\pm\sqrt{m^2-m+1}\in\sigma(D(\overline{\Gamma_{cent}(M_{2mn})})),$ each with multiplicity 1.
\end{enumerate}\et

Therefore, by Theorem~\ref{dm2t1a}, $\overline{\Gamma_{cent}(M_{2mn})}$ is distance integral for even $m,$ if  $\sqrt{m^2-2m+4}$ is even. And it is distance integral for odd $m,$ if $m^2-m+1$ is a perfect square.

Similar to Theorem~\ref{dg2t2}, we have the following.
\bt\label{dlm2t1a} Let $\overline{\Gamma_{cent}(M_{2mn})}$ be the co-centralizer graph of the Metacyclic group of order $2mn$. If $m$ is even, then
\begin{enumerate}
\item[{\rm(a)}] $0\in\sigma(D^L(\overline{\Gamma_{cent}(M_{2mn}))})$ with multiplicity $1;$
\item[{\rm(b)}] $(\frac{m}{2}+1)\in\sigma(D^L(\overline{\Gamma_{cent}(M_{2mn})})),$ with multiplicity $1;$
\item[{\rm(c)}] $(m+1)\in\sigma(D^L(\overline{\Gamma_{cent}(M_{2mn})})),$ with multiplicity $\frac{m}{2}-1;$
\end{enumerate}
And if $m$ is odd, then
\begin{enumerate}
\item[{\rm(a)}] $0\in\sigma(D^L(\overline{\Gamma_{cent}(M_{2mn}))})$ with multiplicity $1;$
\item[{\rm(b)}] $(m+1)\in\sigma(D^L(\overline{\Gamma_{cent}(M_{2mn})})),$ with multiplicity $1.$
\item[{\rm(c)}] $(2m+1)\in\sigma(D^L(\overline{\Gamma_{cent}(M_{2mn})})),$ with multiplicity $m-1.$\end{enumerate}\et

Therefore, $\Gamma_{cent}(M_{2mn})$ is always distance Laplacian integral. Also, similar to Theorem~\ref{dg2t3}, we have the following.

\bt\label{dqm2t1a} Let $\overline{\Gamma_{cent}(M_{2mn})}$ be the co-centralizer graph of the Metacyclic group of order $2mn$. If $m$ is even, then
\begin{enumerate}
\item[{\rm(a)}] $m-3\in\sigma(D^Q(\overline{\Gamma_{cent}(M_{2mn}))})$ with multiplicity $\frac{m}{2}-1;$
\item[{\rm(b)}] $\frac{1}{2}\left[(\frac{5m}{2}-3)\pm\sqrt{\frac{9m^2}{4}-7m+9}\right]\in\sigma(D^Q(\overline{\Gamma_{cent}(M_{2mn})})),$ each with multiplicity $1;$
\end{enumerate}
And if $m$ is odd, then

\begin{enumerate}
\item[{\rm(a)}] $2m-3\in\sigma(D^Q(\overline{\Gamma_{cent}(M_{2mn}))})$ with multiplicity $m-1;$
\item[{\rm(b)}] $\frac{1}{2}\left[(5m-3)\pm\sqrt{9m^2-14m+9}\right]\in\sigma(D^Q(\overline{\Gamma_{cent}(M_{2mn})})),$ each with multiplicity $1.$
\end{enumerate}\et

Therefore, by Theorem~\ref{dqm2t1a}, $\overline{\Gamma_{cent}(M_{2mn})}$ is distance signless Laplacian integral for even $m,$ if  $(\frac{5m}{2}-3)\pm\sqrt{\frac{9m^2}{4}-7m+9}$ is even. And it is distance signless Laplacian integral for odd $m,$ if $(5m-3)\pm\sqrt{9m^2-14m+9}$ is even.

\section{Spectra of $\overline{\Gamma_{cent}(PSL(2,2^k))}$, $k\geq1$}
In this section, we consider the co-centralizer graph of the projective special linear group $PSL(2,2^k)$, where $k\geq2,$  and obtain the distance, distance Laplacian and distance signless Laplacian spectra of it. It follows from \cite{rkn} that $PSL(2,2^k)$ has $2^k+1$ centralizers of cardinality $2^k$, $2^{k-1}(2^k+1)$ centralizers of cardinality $2^k-1$ and $2^{k-1}(2^k-1)$ centralizers of cardinality $2^k+1$. Therefore the co-centralizer graph of $PSL(2,2^k)$ is the complete tripartite graph $K_{2^k+1,2^{k-1}(2^k+1),2^{k-1}(2^k-1)}$.
\subsection{Distance spectrum of $\overline{\Gamma_{cent}(PSL(2,2^k))}$}
By equation (\ref{e1}) we get,
\begin{eqnarray*}
 P_{D(\overline{\Gamma_{cent}(PSL(2,2^k))}}(\lambda)&=&(\lambda+2)^{2^k+2^{2k}-2}[\lambda^3+(4-2^{2k+1}-2^{k+1})\lambda^2+\\
  &&\{4+3\times(2^{4k-2})+3\times(2^{3k})-23\times (2^{2k-2})-2^{k+3}\}\lambda+\\
&&\{-2^{5k}+2^{4k-1}+7\times(2^{3k})-5\times2^{2k-1}-2^{k+3})\}]
\end{eqnarray*}
Therefore, the eigenvalues of $D(\overline{\Gamma_{cent}(PSL(2,2^k)})$ are $-2$ with multiplicity $2^k+2^{2k}-2,$ and three roots of the equation $\lambda^3+(4-2^{2k+1}-2^{k+1})\lambda^2+\{4+3\times(2^{4k-2})+3\times(2^{3k})-23\times (2^{2k-2})-2^{k+3}\}\lambda+\{-2^{5k}+2^{4k-1}+7\times(2^{3k})-5\times2^{2k-1}-2^{k+3})\}=0$.
\subsection{Distance Laplacian spectrum of $\overline{\Gamma_{cent}(PSL(2,2^k))}$}
Let $\1_n$ (resp. $\0_n$) denote the $n\times 1$ vector with each entry 1 (resp. 0). Also, let $J_n$ denote the matrix of order $n$ with all entries equal to 1 (we will write $J$ if the order is clear from the context). It is well known that the distance Laplacian matrix for any graph is positive semidefinite with 0 being the smallest eigenvalue with multiplicity 1 and $\1_n$ is the corresponding eigenvector. The following theorem describes the distance Laplacian spectrum of $\overline{\Gamma_{cent}(PSL(2,2^k))}.$

\bt\label{dlt1} Let $\overline{\Gamma_{cent}(PSL(2,2^k))}$ be the co-centralizer graph of the projective special linear group  $PSL(2,2^k)$. Then
\begin{enumerate}
\item[{\rm(a)}] $0\in\sigma(D^L(\overline{\Gamma_{cent}(PSL(2,2^k))}))$ with multiplicity $1;$
\item[{\rm(b)}] $3\times2^{2k-1}+3\times 2^{k-1}+1\in\sigma(D^L(\overline{\Gamma_{cent}(PSL(2,2^k))})),$ with multiplicity $2^{k-1}(2^k+1)-1;$
\item[{\rm(c)}] $3\times 2^{2k-1}+2^{k-1}+1\in\sigma(D^L(\overline{\Gamma_{cent}(PSL(2,2^k))}))$  with multiplicity $2^{k-1}(2^k-1)-1;$
\item[{\rm(d)}] $2^{k+1}+2^{2k}+2\in\sigma(D^L(\overline{\Gamma_{cent}(PSL(2,2^k))}))$ with multiplicity $2^{k};$
\item[{\rm(e)}] $2^{2k}+2^{k}+1\in\sigma(D^L(\overline{\Gamma_{cent}(PSL(2,2^k))}))$ with multiplicity $2;$
\end{enumerate}\et
 {\bf Proof.} With a suitable labeling of the vertices, the distance Laplacian matrix for $\overline{\Gamma_{cent}(PSL(2,2^k))}$ can be written as
\begin{eqnarray*}
   && D^L(\overline{\Gamma_{cent}(PSL(2,2^k))})= \\
 && \left[
                        \begin{array}{c|c|c}
                          (2^{k+1}+2^{2k}+2)I-2J & -J & -J \\
                          \hline
                          -J & [3(2^{2k-1})+3(2^{k-1})+1]I-2J  & -J \\
                          \hline
                          -J & -J &  [3(2^{2k-1})+2^{k-1}+1]I-2J \\
                        \end{array}
                      \right]
\end{eqnarray*}

Obviously, 0 is an eigenvalue of $D^L(\overline{\Gamma_{cent}(PSL(2,2^k))})$ with multiplicity 1 and $\1_{2^{2k}+2^{k}+1}$ is the corresponding eigenvector.

Also,
\begin{eqnarray*}
  && D^L(\overline{\Gamma_{cent}(PSL(2,2^k))})\left(
                             \begin{array}{c}
                               \0_{2^k+1} \\
                               \hline
                               -1 \\
                               1 \\
                               \0_{2^{k-1}(2^k+1)-2} \\
                               \hline
                               \0_{2^{k-1}(2^k-1)} \\
                             \end{array}
                           \right) \\
  &=& \left(3\times2^{2k-1}+3\times2^{k-1}+1\right)\left(
                             \begin{array}{c}
                                \0_{2^k+1} \\
                               \hline
                               -1 \\
                               1 \\
                               \0_{2^{k-1}(2^k+1)-2} \\
                               \hline
                               \0_{2^{k-1}(2^k-1)} \\
                             \end{array}
                           \right).
\end{eqnarray*}

Therefore, $3\times2^{2k-1}+3\times2^{k-1}+1$ is an eigenvalue of $D^L(\overline{\Gamma_{cent}(PSL(2,2^k))}),$ and in this way we can construct the following set $S_1$ of $2^{k-1}(2^k+1)-1$ independent eigenvectors corresponding to $3\times2^{2k-1}+3\times2^{k-1}+1;$
$$S_1=\left\{\left(
         \begin{array}{c}
           \0_{2^k+1} \\
                               \hline
                               -1 \\
                               1 \\
                               \0_{2^{k-1}(2^k+1)-2} \\
                               \hline
                               \0_{2^{k-1}(2^k-1)} \\
         \end{array}
       \right),\left(
                 \begin{array}{c}
                   \0_{2^k+1} \\
                               \hline
                               -1 \\
                               0 \\
                               1 \\
                               \0_{2^{k-1}(2^k+1)-3} \\
                               \hline
                               \0_{2^{k-1}(2^k-1)} \\
                 \end{array}
               \right),\ldots,\left(
                                \begin{array}{c}
                                  \0_{2^k+1} \\
                               \hline
                               -1 \\
                               \0_{2^{k-1}(2^k+1)-2} \\
                               1 \\
                               \hline
                               \0_{2^{k-1}(2^k-1)} \\
                                \end{array}
                              \right)\right\}.$$
Similarly,
\begin{eqnarray*}
   && D^L(\overline{\Gamma_{cent}(PSL(2,2^k))})\left(
                                   \begin{array}{c}
                                   \0_{2^k+1} \\
                               \hline
                               \0_{2^{k-1}(2^k+1)} \\
                               \hline
                               -1 \\
                               1 \\
                               \0_{2^{k-1}(2^k-1)-2} \\
                                   \end{array}
                                 \right) \\
  &=& \left(3\times2^{2k-1}+2^{k-1}+1\right)\left(
                                   \begin{array}{c}
                                   \0_{2^k+1} \\
                               \hline
                               \0_{2^{k-1}(2^k+1)} \\
                               \hline
                               -1 \\
                               1 \\
                               \0_{2^{k-1}(2^k-1)-2} \\
                                   \end{array}
                                 \right),
\end{eqnarray*}
shows that $3\times2^{2k-1}+2^{k-1}+1$ is an eigenvalue of $D^L(\overline{\Gamma_{cent}(PSL(2,2^k))})$. In this way we can construct the following set $S_2$ of $2^{k-1}(2^k-1)-1$ independent eigenvectors corresponding to $3\times2^{2k-1}+2^{k-1}+1;$
                                           $$S_2=\left\{\left(
                                                   \begin{array}{c}
                                                      \0_{2^k+1} \\
                               \hline
                               \0_{2^{k-1}(2^k+1)} \\
                               \hline
                               -1 \\
                               1 \\
                               \0_{2^{k-1}(2^k-1)-2} \\
                                                   \end{array}
                                                 \right),\left(
                                                           \begin{array}{c}
                                                             \0_{2^k+1} \\
                               \hline
                               \0_{2^{k-1}(2^k+1)} \\
                               \hline
                               -1 \\
                               0 \\
                               1 \\
                               \0_{2^{k-1}(2^k-1)-3} \\
                                                           \end{array}
                                                         \right),\ldots,\left(
                                                                          \begin{array}{c}
                                                                            \0_{2^k+1} \\
                               \hline
                               \0_{2^{k-1}(2^k+1)} \\
                               \hline
                               -1 \\
                               \0_{2^{k-1}(2^k-1)-2} \\
                               1 \\
                                                                          \end{array}
                                                                        \right)\right\}.$$

In a similar way,
\begin{eqnarray*}
  && D^L(\overline{\Gamma_{cent}(PSL(2,2^k))})\left(
                                  \begin{array}{c}
                                    -1 \\
                                     1 \\
                                    \0_{2^k-1} \\
                                    \hline
                                    \0_{2^{k-1}(2^k+1)} \\
                                    \hline
                                    \0_{2^{k-1}(2^k-1)} \\
                                  \end{array}
                                \right) \\
  &=& (2^{k+1}+2^{2k}+2)\left(
                                  \begin{array}{c}
                                     -1 \\
                                     1 \\
                                    \0_{2^k-1} \\
                                    \hline
                                    \0_{2^{k-1}(2^k+1)} \\
                                    \hline
                                    \0_{2^{k-1}(2^k-1)} \\
                                  \end{array}
                                \right).
\end{eqnarray*}

Thus, $2^{k+1}+2^{2k}+2$ is an eigenvalue of $D^L(\overline{\Gamma_{cent}(PSL(2,2^k))}),$ and in this way we can construct the following set $S_3$ of $2^k$ independent eigenvectors corresponding to $2^{k+1}+2^{2k}+2;$
$$S_3=\left\{\left(
        \begin{array}{c}
          -1 \\
          1 \\
          \0_{2^k-1} \\
          \hline
          \0_{2^{k-1}(2^k+1)} \\
          \hline
          \0_{2^{k-1}(2^k-1)}\\
        \end{array}
      \right),\left(
                \begin{array}{c}
                 -1 \\
                 0\\
          1 \\
          \0_{2^k-2} \\
          \hline
          \0_{2^{k-1}(2^k+1)} \\
          \hline
          \0_{2^{k-1}(2^k-1)}\\
                \end{array}
              \right),\ldots,\left(
                               \begin{array}{c}
                                 -1 \\
          \0_{2^k-1} \\
          1 \\
          \hline
          \0_{2^{k-1}(2^k+1)} \\
          \hline
          \0_{2^{k-1}(2^k-1)}\\
                               \end{array}
                             \right)\right\}.$$
Finally, $D^L(\overline{\Gamma_{cent}(PSL(2,2^k))})\left(
                                  \begin{array}{c}
                                    -{2^{k-1}}\1_{2^k+1} \\
                                    \hline
                                    \1_{2^{k-1}(2^k+1)} \\
                                    \hline
                                    \0_{2^{k-1}(2^k-1)} \\
                                  \end{array}
                                \right)=(2^{2k}+2^{k}+1)\left(
                                  \begin{array}{c}
                                     -{2^{k-1}}\1_{2^k+1} \\
                                    \hline
                                    \1_{2^{k-1}(2^k+1)} \\
                                    \hline
                                    \0_{2^{k-1}(2^k-1)} \\
                                  \end{array}
                                \right),$
and
\begin{eqnarray*}
  &&D^L(\overline{\Gamma_{cent}(PSL(2,2^k))})\left(
                                  \begin{array}{c}
                                    -{2^{k-1}(2^k-1)}\1_{2^k+1} \\
                                    \hline
                                    \0_{2^{k-1}(2^k+1)} \\
                                    \hline
                                    (2^k+1)\1_{2^{k-1}(2^k-1)} \\
                                  \end{array}
                                \right) \\
  &=& (2^{2k}+2^{k}+1)\left(
                                  \begin{array}{c}
                                     -2^{k-1}(2^k-1)\1_{2^k+1} \\
                                    \hline
                                    \0_{2^{k-1}(2^k+1)} \\
                                    \hline
                                    (2^k+1)\1_{2^{k-1}(2^k-1)} \\
                                  \end{array}
                                \right).
\end{eqnarray*}

Thus, $2^{2k}+2^{k}+1$ is an eigenvalue of $D^L(\overline{\Gamma_{cent}(PSL(2,2^k))}),$ and in this way we can construct the following set $S_4$ of two independent eigenvectors corresponding to $2^{2k}+2^{k}+1;$
$$S_4=\left\{\left(
        \begin{array}{c}
          {-2^{k-1}}\1_{2^k+1} \\
                                    \hline
                                    \1_{2^{k-1}(2^k+1)} \\
                                    \hline
                                    \0_{2^{k-1}(2^k-1)} \\
        \end{array}
      \right),\left(
                \begin{array}{c}
                 {-2^{k-1}(2^k-1)}\1_{2^k+1} \\
                                    \hline
                                    \0_{2^{k-1}(2^k+1)} \\
                                    \hline
                                    (2^k+1)\1_{2^{k-1}(2^k-1)} \\
                \end{array}
              \right)\right\}.$$

It can be seen that $\1_{2^{2k}+2^{k}+1}\cup S_1\cup S_2\cup S_3\cup S_4$ is a set of mutually orthogonal vectors. Since the order of $\overline{\Gamma_{cent}(PSL(2,2^k))}$ is $2^{2k}+2^{k}+1,$ the result follows.\qed

Hence, from Theorem~\ref{dlt1} it follows that $\overline{\Gamma_{cent}(PSL(2,2^k))}$ is distance Laplacian integral for all $k.$
\subsection{Distance signless Laplacian spectrum of $\overline{\Gamma_{cent}(PSL(2,2^k))}$}

The following theorem describes the distance signless Laplacian spectrum of $\overline{\Gamma_{cent}(PSL(2,2^k))}.$

\bt\label{dhhhqt1} Let $\overline{\Gamma_{cent}(PSL(2,2^k))}$ be the co-centralizer graph of the projective special linear group $PSL(2,2^k)$. Then its distance signless Laplacian spectrum consists of:
\begin{enumerate}
\item[{\rm(a)}] $2^{k+1}+2^{2k}-2$ with multiplicity $2^k,$
\item[{\rm(b)}] $3\times  2^{2k-1}+3\times 2^{k-1}-3$ with multiplicity $2^{k-1}(2^k+1)-1,$
\item[{\rm(c)}] $3\times2^{2k-1}+2^{k-1}+3$  with multiplicity $2^{k-1}(2^k-1)-1,$ and
\item[{\rm (d)}] the three eigenvalues of the matrix $$\mathfrak{D}_P=\tiny{\left[
                                                        \begin{array}{c|c|c}
                                                          2^{k+1}+2^{2k}-2+2(2^k+1) & 2^{k-1}(2^k+1) & 2^{k-1}(2^k-1) \\
                                                          \hline
                                                          2^k+1 & 3\times2^{2k-1}+3\times 2^{k-1}-3+2^k(2^k+1)& 2^{k-1}(2^k-1) \\
                                                          \hline
                                                          2^k+1 & 2^{k-1}(2^k+1) &  3\times 2^{2k-1}+2^{k-1}-3 +2^k(2^k-1) \\
                                                        \end{array}
                                                      \right]}.
$$
\end{enumerate}\et
{\bf Proof.} With a suitable labeling of the vertices, the distance signless Laplacian matrix for $\overline{\Gamma_{cent}(PSL(2,2^k))}$ can be written as
$$D^Q(\overline{\Gamma_{cent}(PSL(2,2^k))})=$$
$$\left[
                        \begin{array}{c|c|c}
                          (2^{k+1}+2^{2k}-2)I+2J & J & J \\
                          \hline
                          J & (3(2^{2k-1})+3(2^{k-1})-3)I+2J  & J \\
                          \hline
                          J & J & (3(2^{2k-1})+2^{k-1}-3)I+2J \\
                        \end{array}
                      \right].$$

    Now, $D^Q(\overline{\Gamma_{cent}(PSL(2,2^k))})\left(
                                  \begin{array}{c}
                                    -1 \\
                                     1 \\
                                    \0_{2^k-1} \\
                                    \hline
                                    \0_{2^{k-1}(2^k+1)} \\
                                    \hline
                                    \0_{2^{k-1}(2^k-1)} \\
                                  \end{array}
                                \right)=(2^{k+1}+2^{2k}-2)\left(
                                  \begin{array}{c}
                                     -1 \\
                                     1 \\
                                    \0_{2^k-1} \\
                                    \hline
                                    \0_{2^{k-1}(2^k+1)} \\
                                    \hline
                                    \0_{2^{k-1}(2^k-1)} \\
                                  \end{array}
                                \right).$

Thus, $2^{k+1}+2^{2k}-2$ is an eigenvalue of $D^Q(\overline{\Gamma_{cent}(PSL(2,2^k))}),$ and in this way we can construct the following set $S_1$ of $2^k$ independent eigenvectors corresponding to $2^{k+1}+2^{2k}-2;$
$$S_1=\left\{\left(
        \begin{array}{c}
          -1 \\
          1 \\
          \0_{2^k-1} \\
          \hline
          \0_{2^{k-1}(2^k+1)} \\
          \hline
          \0_{2^{k-1}(2^k-1)}\\
        \end{array}
      \right),\left(
                \begin{array}{c}
                 -1 \\
                 0\\
          1 \\
          \0_{2^k-2} \\
          \hline
          \0_{2^{k-1}(2^k+1)} \\
          \hline
          \0_{2^{k-1}(2^k-1)}\\
                \end{array}
              \right),\ldots,\left(
                               \begin{array}{c}
                                 -1 \\
          \0_{2^k-1} \\
          1 \\
          \hline
          \0_{2^{k-1}(2^k+1)} \\
          \hline
          \0_{2^{k-1}(2^k-1)}\\
                               \end{array}
                             \right)\right\}.$$

Also,
\begin{eqnarray*}
  &&D^Q(\overline{\Gamma_{cent}(PSL(2,2^k))})\left(
                             \begin{array}{c}
                               \0_{2^k+1} \\
                               \hline
                               -1 \\
                               1 \\
                               \0_{2^{k-1}(2^k+1)-2} \\
                               \hline
                               \0_{2^{k-1}(2^k-1)} \\
                             \end{array}
                           \right) \\
  &=& \left(3\times 2^{2k-1}+3\times2^{k-1}-3\right)\left(
                             \begin{array}{c}
                                \0_{2^k+1} \\
                               \hline
                               -1 \\
                               1 \\
                               \0_{2^{k-1}(2^k+1)-2} \\
                               \hline
                               \0_{2^{k-1}(2^k-1)} \\
                             \end{array}
                           \right).
\end{eqnarray*}
Therefore, $3\times 2^{2k-1}+3\times2^{k-1}-3$ is an eigenvalue of $D^Q(\overline{\Gamma_{cent}(PSL(2,2^k))}),$ and in this way we can construct the following set $S_2$ of $2^{k-1}(2^k+1)-1$ independent eigenvectors corresponding to $3\times 2^{2k-1}+3\times2^{k-1}-3;$
$$S_2=\left\{\left(
         \begin{array}{c}
           \0_{2^k+1} \\
                               \hline
                               -1 \\
                               1 \\
                               \0_{2^{k-1}(2^k+1)-2} \\
                               \hline
                               \0_{2^{k-1}(2^k-1)} \\
         \end{array}
       \right),\left(
                 \begin{array}{c}
                   \0_{2^k+1} \\
                               \hline
                               -1 \\
                               0 \\
                               1 \\
                               \0_{2^{k-1}(2^k+1)-3} \\
                               \hline
                               \0_{2^{k-1}(2^k-1)} \\
                 \end{array}
               \right),\ldots,\left(
                                \begin{array}{c}
                                  \0_{2^k+1} \\
                               \hline
                               -1 \\
                               \0_{2^{k-1}(2^k+1)-2} \\
                               1 \\
                               \hline
                               \0_{2^{k-1}(2^k-1)} \\
                                \end{array}
                              \right)\right\}.$$

Similarly,
\begin{eqnarray*}
  &&D^Q(\overline{\Gamma_{cent}(PSL(2,2^k))})\left(
                                   \begin{array}{c}
                                   \0_{2^k+1} \\
                               \hline
                               \0_{2^{k-1}(2^k+1)} \\
                               \hline
                               -1 \\
                               1 \\
                               \0_{2^{k-1}(2^k-1)-2} \\
                                   \end{array}
                                 \right) \\
  &=& \left(3\times2^{2k-1}+2^{k-1}+3\right)\left(
                                   \begin{array}{c}
                                   \0_{2^k+1} \\
                               \hline
                               \0_{2^{k-1}(2^k+1)} \\
                               \hline
                               -1 \\
                               1 \\
                               \0_{2^{k-1}(2^k-1)-2} \\
                                   \end{array}
                                 \right),
\end{eqnarray*}
shows that $3\times2^{2k-1}+2^{k-1}+3$ is an eigenvalue of $D^Q(\overline{\Gamma_{cent}(PSL(2,2^k))})$. In this way we can construct the following set $S_3$ of $2^{k-1}(2^k-1)-1$ independent eigenvectors corresponding to $3\times2^{2k-1}+2^{k-1}+3;$
                                           $$S_3=\left\{\left(
                                                   \begin{array}{c}
                                                      \0_{2^k+1} \\
                               \hline
                               \0_{2^{k-1}(2^k+1)} \\
                               \hline
                               -1 \\
                               1 \\
                               \0_{2^{k-1}(2^k-1)-2} \\
                                                   \end{array}
                                                 \right),\left(
                                                           \begin{array}{c}
                                                             \0_{2^k+1} \\
                               \hline
                               \0_{2^{k-1}(2^k+1)} \\
                               \hline
                               -1 \\
                               0 \\
                               1 \\
                               \0_{2^{k-1}(2^k-1)-3} \\
                                                           \end{array}
                                                         \right),\ldots,\left(
                                                                          \begin{array}{c}
                                                                            \0_{2^k+1} \\
                               \hline
                               \0_{2^{k-1}(2^k+1)} \\
                               \hline
                               -1 \\
                               \0_{2^{k-1}(2^k-1)-2} \\
                               1 \\
                                 \end{array}
                                \right)\right\}.$$
Thus, we have obtained $2^k+2^{k-1}(2^k+1)-1+2^{k-1}(2^k-1)-1=2^k+2^{2k}-2$ eigenvalues of $D^Q(\overline{\Gamma_{cent}(PSL(2,2^k))}).$ Moreover, we note that all the eigenvectors constructed so far,
are orthogonal to $\left[
                     \begin{array}{c}
                       \1_{2^k+1} \\
                       \hline
                       \0_{2^{k-1}(2^k+1)} \\
                       \hline
                       \0_{2^{k-1}(2^k-1)}\\
                     \end{array}
                   \right],~\left[
                     \begin{array}{c}
                       \0_{2^k+1} \\
                       \hline
                       \1_{2^{k-1}(2^k+1)} \\
                       \hline
                       \0_{2^{k-1}(2^k-1)}\\
                     \end{array}
                   \right]$ and $\left[
                     \begin{array}{c}
                       \0_{2^k+1} \\
                       \hline
                       \0_{2^{k-1}(2^k+1)} \\
                       \hline
                       \1_{2^{k-1}(2^k-1)}\\
                     \end{array}
                   \right].$
Therefore, these three vectors span the remaining three eigenvectors of $D^Q(\overline{\Gamma_{cent}(PSL(2,2^k))}).$ Thus, the remaining eigenvectors
of $D^Q(\overline{\Gamma_{cent}(PSL(2,2^k))})$ are of the form $\left[
                     \begin{array}{c}
                       a\1_{2^k+1} \\
                       \hline
                       b\1_{2^{k-1}(2^k+1)} \\
                       \hline
                       c\1_{2^{k-1}(2^k-1)}\\
                     \end{array}
                   \right],$
for some $(a,b,c)\ne(0,0,0).$ Therefore if $\mu$ is an eigenvalue of $D^Q(\overline{\Gamma_{cent}(PSL(2,2^k))})$ with eigenvector $\left[
                     \begin{array}{c}
                       a\1_{2^k+1} \\
                       \hline
                       b\1_{2^{k-1}(2^k+1)} \\
                       \hline
                       c\1_{2^{k-1}(2^k-1)}\\
                     \end{array}
                   \right],$
then $a,b,c$ are the solution of the following system of equation
\begin{eqnarray*}
  \left(2^{k+1}+2^{2k}-2+2\times(2^k+1)\right)a+\left(2^{k-1}\times (2^k+1)\right)b+\left(2^{k-1}\times(2^k-1)\right)c &=& 0 \\
   (2^k+1)a +\left(3\times2^{2k-1}+3\times2^{k-1}-3+2\times2^{k-1}\times(2^k+1)\right)b+\left(2^{k-1}\times(2^k-1)\right)c&=& 0 \\
    (2^k+1)a+\left(2^{k-1}\times(2^k+1)\right)b+\left(3\times2^{2k-1}+2^{k-1}-3+2\times2^{k-1}(2^k-1)\right)c&=& 0.
\end{eqnarray*}
Therefore, the remaining three eigenvalues of $D^Q(\overline{\Gamma_{cent}(PSL(2,2^k))})$ are the eigenvalues of the matrix
$\mathfrak{D}_P.$\qed

Hence, by Theorem~\ref{dhhhqt1}, $D^Q(\overline{\Gamma_{cent}(PSL(2,2^k))})$ is distance signless Laplacian integral if $\mathfrak{D}_P$ have integral spectrum.

\section{Conclusion}
In this article, we have investigated the distance, distance (signless) Laplacian spectra of co-centralizer graphs of the generalized quaternion group, the dihedral group, the quasidihedral group, the metacyclic group, and the projective special linear group. We also obtain conditions under which these graphs will be distance, distance (signless) Laplacian integral.
\end{document}